\documentclass{emsprocart}

\contact[Christoph.Schweigert@uni-hamburg.de]
{Algebra und Zahlentheorie, Fachbereich Mathematik, 
Universit\"at Hamburg, Bundesstra\ss{}e 55, 
20148 Hamburg, Germany}

\contact[juergen.fuchs@kau.se ]
{Teoretisk fysik, Karlstads universitet, Universitetsgatan 21,
651\,88 Karlstad, Sweden}

\newtheorem{theorem}{Theorem}[section]

\newtheorem{proposition}[theorem]{Proposition}

\theoremstyle{definition}
\newtheorem{definition}[theorem]{Definition}

\title[Low-dimensional geometry and representation theory]
{Low-dimensional topology, low-dimensional field \\
theory and representation theory}

\author[J.\ Fuchs and C.\ Schweigert]
{J\"urgen Fuchs and Christoph Schweigert
\thanks{CS is partially supported by the Collaborative Research Centre 676 ``Particles, 
Strings and the Early Universe - the Structure of Matter and Space-Time'', by 
the RTG 1670 ``Mathematics inspired by String theory and Quantum Field Theory'' 
and by the DFG Priority Programme 1388 ``Representation Theory''.
\\
JF is supported by VR under project no.\ 621-2013-4207.
\\
We thank Simon Lentner and Chris Schommer-Pries for helpful comments on the manuscript.
}}

\begin{document}

\begin{abstract}
Structures in low-dimensional topology and low-dimensional geometry -- often 
combined with ideas from (quantum) field theory --  can explain and inspire
concepts in algebra and in representation theory and their categorified versions. 
We present a personal view on some of these instances which have appeared within the
Research Priority Programme SPP 1388 ``Representation theory''.
\end{abstract}



\maketitle


\section{Introduction}

Structures and relations in algebra and representation theory are sometimes
``explained'' by geometric or topological facts, frequently by facts concerning
a category of low-dimensional geometric or topological objects. 

As a first illustration, consider the following well-known example of an algebraic relation.
Let $A$ be a set. A map $m\colon A\,{\times}\, A \,{\to}\, A$ 
is called {\em associative} if for any integer $n \,{\geq}\, 3$ and any given $n$-tuple
$(a_1,a_2\ldots,a_n) \,{\in}\, A^n$, every way of putting brackets in the expression
$a_1\, a_2\, \cdots \,a_n$ gives, upon applying $m$ repeatedly according to the chosen
bracketing, one and the same element of $A$. This is not the typical
the textbook definition of associativity; the latter is rather based on the

\begin{proposition}
A map $m\colon A\,{\times}\, A \,{\to}\, A$ is associative if and only if the equality
  $$
  m(m(a_1, a_2), a_3) = m( a_1,m (a_2, a_3))
  $$
holds for all triples $a_1,a_2,a_3 \,{\in}\, A$.
\end{proposition}

It should be appreciated that this is really an assertion about binary trees with 
a finite number of leaves. Indeed, any bracketing of $n$ elements corresponds to a 
binary tree with $n$ leaves. The statement then follows from the observation that any two
binary trees can be transformed into each other by applying repeatedly a move 
that amounts to the associativity relation. 
In this sense, a combinatorial  property of the collection 
of all binary trees ``explains'' the textbook definition of associativity.

To give an example of a structure ``explained'' or, rather, justified, by facts from 
geometry and related to quantum field theory, consider the symmetric monoidal category
$\mathcal C\!ob(2,1)$ of smooth oriented cobordisms. The objects of $\mathcal C\!ob(2,1)$ 
are finite disjoint unions of oriented circles ${\mathbb S}^1$, and its morphisms are 
diffeomorphism classes of smooth oriented surfaces with boundary.
One way to learn more about the category $\mathcal C\!ob(2,1)$ is to {\em represent} it
on a symmetric monoidal category $\mathcal S$, i.e.\ to study symmetric monoidal functors
$\mathrm{tft}\colon \mathcal C\!ob(2,1) \,{\to}\, {\mathcal S}$. (An important example 
is ${\mathcal S} \,{=}\, \mathrm{vect}$, the category of
vector spaces over a field $\mathbb K$.) Such a representation $\mathrm{tft}$ of 
the category $\mathcal C\!ob(2,1)$ is called a (2,1)-dimensional topological field theory.

For the rest of this section we restrict our attention to $\mathbb K$-linear
representations, with $\mathbb K$ a field, i.e\ to 
the case of ${\mathcal S} \,{=}\, \mathrm{vect}$.
Representations can be constructed using generators and relations for the object to be 
represented. The category $\mathcal C\!ob(2,1)$ admits such a description in
terms of generators and relations, given by the pair-of-pants decomposition of 
Riemann surfaces with boundary. As a consequence, it suffices to know the vector space 
  \begin{equation}\label{tftS1=C}
  \mathrm{tft}({\mathbb S}^1) =: C
  \end{equation}
associated to the circle together with the linear maps $m\colon C \,{\otimes}\, C \,{\to}\, C$ 
and $\Delta\colon C \,{\to}
    $\linebreak[0]$
C\,{\otimes}\, C$ that are obtained from the three-punc\-tured sphere 
regarded as a cobordism 
${\mathbb S}^1 \,{\sqcup}\, {\mathbb S}^1 \,{\to}\, {\mathbb S}^1$ and as
${\mathbb S}^1\,{\to}\, {\mathbb S}^1 \,{\sqcup}\, {\mathbb S}^1$, respectively, as well as the
maps $\eta\colon {\mathbb K}\,{\to}\, C$ and $\varepsilon\colon C \,{\to}\, {\mathbb K}$ 
obtained from the disk, seen as a cobordism $\emptyset \,{\to}\, {\mathbb S}^1$ and as
${\mathbb S}^1 \,{\to}\, \emptyset$, respectively.

It is well known (see e.g.\ the textbook \cite{cs-Kock}) that the relations in the 
pair-of-pants decompositions imply that $(C,m,\Delta,\eta,\epsilon)$ is a commutative 
Frobenius algebra over $\mathbb K$. In this sense, the structure of $\mathcal C\!ob(2,1)$
together with the idea, inspired from quantum field theory, to study representations of 
$\mathcal C\!ob(2,1)$, enforce on us the algebraic notion of a commutative Frobenius algebra 
(which was, of course, a known structure long before the advent of topological
field theory). The result amounts to a classification of
topological field theories, which are objects in the functor category
$ F\!un_{\otimes,\mathrm{sym}}(\mathcal C\!ob(2,1),\mathrm{vect})$: 
         
\begin{proposition}
The evaluation on the circle provides an equivalence
  \begin{equation}
  F\!un_{\otimes,\mathrm{sym}}(\mathcal C\!ob(2,1),\mathrm{vect}) \,\simeq\, F\!rob 
  \label{Frobequiv}
  \end{equation}
of categories, where $F\!rob$ is the groupoid of commutative Frobenius algebras
over the field $\mathbb K$.
\end{proposition}

Non-commutative Frobenius algebras and their representation categories can be obtained 
in a similar spirit when one allows for a larger geometric category 
$\mathcal C\!ob(2,1)^\partial$ that has as objects disjoint unions of circles as well
as intervals. One then obtains a so-called open/closed topological
field theory. For details, we refer to \cite{cs-LP}.

\medskip

In the sequel we will present a small, and strongly biased by personal taste, digest of
instances in which the interplay of geometric and algebraic structures has lead to a
mathematical insight. All those instances are, in some way, related to
the Priority Programme ``Representation Theory''.


\section{Radford's $S^4$-theorem and the orthogonal group $\mathrm{SO}(3)$}

Next we present a more recent example for how topological field theory can ``explain''
a classical fact from algebra. The following statements
about a fini\-te-di\-men\-sional Hopf algebra $H$ over a field $\mathbb K$ can be found 
in any textbook on Hopf algebras. The Hopf algebra $H$ has a one-dimensional subspace 
  $$
  I_l := \{ t \,{\in}\, H \,|\, h \,{\cdot}\, t \,{=}\, \varepsilon(h)\,t \}
  $$
of left integrals, with $\varepsilon\colon H \,{\to}\, \mathbb{K}$ the counit of $H$. 
By associativity, for $t \,{\in}\, I_l$ and any $h \,{\in}\, H$,
the element $t \,{\cdot}\, h \,{\in}\, H$ is again a left integral of $H$.
Since the space of integrals is one-dimensional, this in turn implies the existence of  
a linear form $\alpha \,{\in}\, H^*$ satisfying $t \,{\cdot}\, h \,{=}\, \alpha(h)\, t$ 
for any $t \,{\in}\, I_l(H)$. This linear form is a morphism of algebras and thus a
group-like element of $H^*$. By applying the same reasoning to the dual Hopf algebra 
$H^*$ one finds analogously a distinguished group-like element $a \,{\in}\, H$. A classical result 
by Radford is

\begin{theorem}[Radford \cite{cs-rad}]
Let $H$ be a finite-dimensional Hopf algebra over a field $\mathbb K$, and let 
$a \,{\in}\, H$ and $\alpha \,{\in}\, H^*$ be the distinguished group-like elements. 
Then the fourth power of the antipode $S$ of $H$ satisfies
  $$
  S^4(h) = a\, (\alpha^{-1}{\rightharpoonup}\, h\,{\leftharpoonup}\,\alpha)\, a^{-1}
  = \alpha^{-1} \,{\rightharpoonup}\, (aha^{-1}) \,{\leftharpoonup}\, \alpha
  $$
for all $h \,{\in}\,H$.
\end{theorem}

\noindent
Here the symbol $\rightharpoonup$ denotes the left action of $H^*$ on $H$ which, using the
Sweedler notation $\Delta(h) \,{=}\, h_{(1)} \,{\otimes}\, h_{(2)}$ for the coproduct 
of $H$, is given by 
$\alpha\,{\rightharpoonup}\, h \,{=}\, h_{(1)}\, \langle\alpha,h_{(2)}\rangle$,
and analogously $\leftharpoonup$ is the right $H^*$-action on $H$.
There is a purely algebraic proof for the purely algebraic statement in Theorem 2.1. 

The statement admits a categorical reformulation. Consider a finite tensor category 
$\mathcal C$, i.e.\ an abelian rigid tensor category over a field $\mathbb K$ in 
which the monoidal unit is simple and which obeys several finiteness conditions 
(finite-dimensional morphism spaces, finite Jordan-H\"older length for any
object, and finitely many isomorphism classes of simple objects each of 
which has a projective cover). As an example, take the rigid monoidal category
${\mathcal C} \,{=}\, H$-mod of finite-dimensional modules over
a finite-dimensional Hopf algebra $H$. The group-like element $\alpha \,{\in} H^*$ 
is a morphism of algebras and thus defines an invertible object $D$ in 
${\mathcal C} \,{=}\, H$-mod. Such an invertible object can be 
introduced for any finite tensor category \cite[Lemma\,2.9]{cs-etos}. Taking into account 
that the left dual of a module $M \,{=}\, (M,\rho) \,{\in}\, H$-mod
is given by ${}^{\vee}\! M \,{=}\, (M^{*},\rho\,{\circ}\, S^{-1})$ and the right dual 
of $M$ by $M^\vee \,{=}\, (M^{*},\rho\,{\circ}\, S)$, Radford's theorem amounts
to a relation between the double duals ${}^{\vee\vee}\!M$ and $M^{\vee\vee}$ 
in ${\mathcal C} \,{=}\, H$-mod. Indeed, the assertion in Theorem 2.1 
generalizes
\cite{cs-etno2} to the following purely categorical statement:

\begin{theorem}
Let $\mathcal C$ be a finite tensor category. Then there is a natural isomorphism 
  $$ 
  {}^{\vee\vee} ? \longrightarrow D^{-1}\,\otimes\, ?^{\vee\vee\!}\otimes D
  $$
of monoidal functors.
\end{theorem}

A proof of this statement can be given entirely in a category theoretic setting, see
\cite{cs-etno2} and \cite[Sect.\,5.2]{cs-shimi7}.

We now outline how the theorem can be ``understood'' in the framework of topological field
theories. 
The definition of a topological field theory as a
symmetric monoidal functor can be extended to other cobordism categories,
and to higher-ca\-te\-gorical versions of cobordism categories. In a rather general 
framework one can consider an $(\infty,n)$-category $\mathcal C\!ob_{\infty,n}^{\mathrm{fr}}$ 
of $n$-dimensional framed cobordisms, with objects being $n$-framed points, 1-morphisms
$n$-framed 1-manifolds with boundary, 2-morphisms $n$-framed 2-manifolds with corners, 
etc. For $k \,{>}\, n$, all $k$-morphisms are invertible.

A topological field theory with values in $\mathcal S$ is then a symmetric monoidal functor 
$\mathrm{tft}\colon \mathcal C\!ob_{\infty,n}^{\mathrm{f}r}\,{\to}\, {\mathcal S}$, with
$\mathcal S$ a symmetric monoidal $(\infty,n)$-category. The same reasoning that shows 
that the vector space \eqref{tftS1=C} carries additional structure can be used to see that
$\mathrm{tft}(*)$, the object assigned to a point, has special properties.
The cobordism theorem \cite{cs-Lurie} states that $\mathrm{tft}(*)$
is a fully dualizable object in 
$\mathcal S$ and that any such object determines a framed topological field theory. 
There is an equivalence 
  \begin{equation}
  F\!un(\mathcal C\!ob_{\infty,n}^{\mathrm{fr}},\mathcal{S})
  \,\simeq\,  k[\mathcal{S}^{\mathrm{f.d.}}] 
  \end{equation}
of $\infty$-groupoids generalizing (\ref{Frobequiv}), where on the right hand side 
the symbol $k$ indicates that one discards all non-invertible morphisms in 
the full subcategory of fully dualizable objects of $\mathcal S$.

We now focus on three-dimensional topological field theories. An important target
$\mathcal S$ is then the symmetric monoidal 3-category $\mathcal Bimod$
which has finite tensor categories as objects, 
bimodule categories over finite tensor categories as 1-morphisms, and
bimodule functors and bimodule natural transformations as 2- and 3-morphisms,
respectively. A fusion category, i.e.\ a semisimple finite tensor category, is a 
fully dualizable object in $\mathcal Bimod$ and thus determines a
framed three-di\-men\-sional topological field theory. Now by change of framing, 
the orthogonal group $\mathrm{O}(3)$ acts on $\mathcal C\!ob_{\infty,3}^{\mathrm{fr}}$ 
and thus on the functor category $F\!un(\mathcal C\!ob_{\infty,3}^{\mathrm{fr}},
\mathcal Bimod)$.
This translates into a homotopy action of $\mathrm{O}(3)$ on the fully dualizable objects,
i.e.\ on fusion categories: points in $\mathrm{O}(3)$ give self-equivalences 
$k[\mathcal Bimod^{\mathrm{f.d.}}] \,{\to}\, k[\mathcal Bimod^{\mathrm{f.d.}}]$,
paths in $\mathrm{O}(3)$ give natural transformations between self-equivalences, etc. 
The homotopy groups of the Lie group $\mathrm{O}(3)$ are well-known:
The group of connected components is $\pi_0(\mathrm{O}(3)) \,{=}\, \mathbb{Z}_2$; 
the non-tri\-vial component acts on monoidal categories as
$(\mathcal{C},\otimes) \,{\mapsto}\, (\mathcal{C},\otimes^{\mathrm{opp}})$. 
The fundamental group is $\pi_1(\mathrm{O}(3)) \,{=}\, \mathbb{Z}_2$;
its non-tri\-vial element acts on fusion categories as an autoequivalence,
which turns out \cite{cs-DSS} to be given by the bidual $?^{\vee\vee}$. 
Thus from the group-theoretical fact that the non-trivial element of
$\pi_1(\mathrm{O}(3)) \,{=}\, \mathbb{Z}_2$ has order 2 one concludes that 
the quadruple dual $?^{\vee\vee\vee\vee}$ is trivial, which in turn can
be seen to imply Radford's theorem for fusion categories. (Further, there are 
weakenings on the side of topological field theories which allow one to make statements
about the quadruple dual of finite tensor categories, see \cite{cs-DSS}.) In this way, 
topological field theory provides a highly surprising connection between homotopy groups
of Lie groups and Radford's $S^4$-theorem for finite-dimensional Hopf algebras.


\section{Modularization via equivariant Dijkgraaf-Witten theories}

In the previous section we used topological field theories in three dimensions
that are fully extended, i.e.\ down to the point. In the sequel we study a different variant of 
three-di\-men\-sional topological field theories, which is only extended down to
one-di\-men\-sional manifolds. It is of independent interest for representation theory.

Consider the symmetric monoidal bicategory $\mathcal C\!ob(3,2,1)$
whose objects are closed oriented smooth one-dimensional manifolds, 1-morphisms are
oriented surfaces with boundary and 2-morphisms are oriented three-manifolds with corners.
To represent this bicategory, consider a symmetric monoidal 2-functor
$\mathrm{tft}\colon \mathcal C\!ob(3,2,1) 
    $\linebreak[0]$
{\to}\, {\mathcal S}$ with values in 
a symmetric monoidal bicategory $\mathcal S$: this is a called a 3-2-1-ex\-tended 
topological field theory. One possible, and particularly important, choice for 
$\mathcal S$ is the symmetric monoidal bicategory $2$-$\mathrm{vect}$ of
so-called 2-vector spaces, which has finite semisimple $\mathbb K$-linear abelian 
categories as objects, $\mathbb K$-linear functors as 1-morphisms
and natural transformations as 2-morphisms.

In analogy to the situation with $\mathcal C\!ob(2,1)$, the category 
$\mathrm{tft}({\mathbb S}^1)$ is endowed with additional algebraic structure, and
that structure should determine the three-di\-men\-sional topological field theory.
Analyzing this structure on $\mathrm{tft}({\mathbb S}^1)$ 
with the help of a presentation of $\mathcal C\!ob(3,2,1)$ 
in terms of generators and relations \cite{cs-BDSV}, the category 
$\mathrm{tft}({\mathbb S}^1)$ turns out to be a modular tensor category, i.e.\ 
a braided monoidal semisimple category, together with dualities and a ribbon twist,
with non-degenerate braiding. The non-degeneracy condition on the braiding can be 
formulated in several equivalent ways. A converse statement, that a modular tensor 
category determines a 3-2-1-ex\-tended topological field theory,
has already been established long ago \cite{cs-retu}.

Topological field theory thus suggests to regard modular tensor categories 
as a categorified notion of (semisimple) commutative Frobenius algebras. This
is indeed a fruitful point of view.

There are several representation-theoretic sources of modular tensor categories:
\\[-1.5em]
\begin{itemize}\addtolength\itemsep{-5pt}
\item 
Left modules over a connected factorizable ribbon weak Hopf algebra
with Haar integral over an algebraically closed field \cite{cs-NTV}.
\item 
Local sectors of a net of von Neumann algebras on $\mathbb R$ that has finite $\mu$-index 
and is strongly additive and split \cite{cs-KLM}.
\item 
Representations of a selfdual $C_2$-cofinite vertex algebra
with an additional finiteness condition on the homogeneous
components and with semisimple representation category \cite{cs-H}.
\end{itemize}
~\\[-1.6em]
(The last two items constitute  different mathematical formalizations
of chiral conformal field theories.)
Despite these many sources, modular tensor categories are ``rare'' 
objects.\footnote{~This statement should not be interpreted in the sense that 
  modular tensor categories can be classified in a naive way. 
  Indeed, the class of such categories contains the representation categories of all 
  doubles of finite groups, so that a naive classification is impossible. Still, imposing
  a natural equivalence relation, which amounts to dividing out Drinfeld doubles, one 
  arrives at a Witt group of modular tensor categories \cite{cs-dmno} which has a strong
  arithmetic flavor. In applications, this Witt group is a useful
  recipient of obstructions, see \cite{cs-FSV1}.} 
As a potential further source of modular tensor categories,
crossed modules of finite groups have been proposed \cite{cs-ba} some time ago.
A crossed module consists of two finite groups $G_1$ and $G_2$, an action of
$G_2$ on $G_1$ by group automorphisms, and a group homomorphism 
$\partial\colon G_1 \,{\to}\, G_2$ such that
  $$
  \partial(g.m)= g\cdot \partial(m)\cdot g^{-1} \qquad\text{and}\qquad
  (\partial n).m=n^{-1}\cdot m\cdot n
  $$
for all $g \,{\in}\, G_2$ and $m,n \,{\in} G_1$. These data determine a braided monoidal 
category $\mathcal{C}(G_1,G_2,\partial)$ whose objects are $G_1$-graded finite-dimensional 
vector spaces carrying an action of $G_2$ such that
$g(V_m) \,{\subseteq}\, V_{g.m}$, with the braiding given by
  $$
  v_m\otimes v_n \,\longmapsto\, \partial(m).v_n \otimes v_m
  $$
for $v_n$ and $v_m$ homogeneous elements of degree $n,m \,{\in}\, G_1$, respectively. The 
category $\mathcal{C}(G_1,G_2,\partial)$ is a premodular tensor category, i.e.\ has all 
the properties of a modular tensor category except that the braiding may be degenerate.
It is modular if and only if $\partial$ is an isomorphism, in which case
one obtains the representation category of the double of a finite group,
i.e.\ an already known type of modular tensor category. In the general case
the category can be modularized according to a standard procedure \cite{cs-brug,cs-mug}. 
Somewhat disappointingly, also this modularization turns out to be braided
equivalent to the representation category of the Drinfeld double of a finite group.

The modularization leads to a category with a certain weak group action,
suggesting that it is actually the neutral part of an equivariant modular category
\cite{cs-Turaev}. To find equivariant extensions of a given modular tensor category 
$\mathcal C$ is an interesting algebraic problem with deep links to the Brauer-Picard group 
\cite{cs-enoM} of $\mathcal C$. Again topological field theory, in this case a twisted 
version of Dijkgraaf-Witten theories, provides one solution to this algebraic problem.

As a simplifying assumption, let us suppose that the group homomorphism
$\partial\colon G_1 \,{\to}\, G_2$ appearing in the
crossed module of finite groups is an injective group homomorphism, so that the group
$G_1$ is a normal subgroup of $G_2$. The corresponding sequence 
$1 \,{\to}\, G_1 \,{\to}\, G_2 \,{\stackrel\pi\to}\, J \,{\to}\, 1$ of groups is 
not necessarily split. But still the group $J$ acts weakly on $G_1$, which can be seen 
as follows. Choose a set-theoretic section $s \colon J \,{\to}\, G_2$ of $\pi$. Then for
any $j \,{\in}\, J$ the automorphism $\rho_j(g_1) \,{:=}\, s(j)\, g_1 \,s(j)^{-1}$ of $G_1$ 
obeys $\rho_{j_1} {\circ}\, \rho_{j_2} \,{=}\, \mathrm{Ad}_{c_{j_1,j_2}} {\circ}\,
\rho_{j_1\cdot j_2}$ with $c_{ij} \,{:=}\, s(i) \,{\cdot}\, s(j) \,{\cdot}\,
s(i{\cdot}j)^{-1} {\in}\, G_1$. (We refrain from exhibiting the
coherence relations satisfied by these group elements $c_{ij}$.)
This weak action induces a (weak) $J$-action on the modularization of the category
$\mathcal{C}(G_1,G_2,\partial)$ which in the case at hand is the category 
$\mathcal{D}(G_1)$-mod of modules over the double of $G_1$.

For any finite group $G$, the Dijkgraaf-Witten theory provides 
a 3-2-1-extended topological field theory. The construction is based on the following
facts (for a recent presentation, see \cite{cs-morton} and references therein). 
To any manifold $X$ there is associated the groupoid $\mathcal{A}_G(X)$ of $G$-principal
bundles. This groupoid is essentially finite if $X$ is compact. (Extended) cobordisms 
are spans of manifolds, possibly with corners; pull back of bundles yields spans of 
groupoids. Taking functor categories into the symmetric monoidal category $\mathrm{vect}$ 
amounts to a linearization
and results in a 3-2-1-extended topological field theory. In particular one obtains
  $$
  \mathrm{tft}({\mathbb S}^1) = F\!un(\mathcal{A}_G(\mathbb{S}^1),
  \mathrm{vect}) \,\simeq\, \mathcal{D}(G)\mbox{-mod} \, . $$

The main idea of the construction of a $J$-equivariant modular category
whose neutral component is the modularization of $\mathcal{C}(G_1,G_2,\partial)$
is geometric \cite{cs-mns}:
replace the categories $\mathcal{A}_G(X)$ by categories of twisted bundles.

\begin{definition}\label{def:twisted}
Let $J$ act weakly on $G_1$, such that there is an exact sequence
$1 \,{\to}\, G_1 \,{\to}\, G_2 \,{\stackrel\pi\to}\, J \,{\to}\, 1$ of finite groups,
and let $P \,{\stackrel J\to}\, M$ be a $J$-cover over a smooth manifold $M$.
\\[3pt]
(i)\,
A \emph{$P$-twisted $G_1$-bundle over $M$} is a pair $(Q,\varphi)$ consisting of a 
$G_2$-bundle 
$Q$ over $M$ and a smooth map $\varphi\colon Q \,{\to}\, P$ over $M$ that obeys
  $$
  \varphi(q\,{\cdot}\, h) = \varphi(q) \cdot \pi(h) \,\,
  $$
for all $q \,{\in}\, Q$ and all $h \,{\in}\, G_2$. 
\\[3pt]
(ii)\,
A \emph{morphism $(Q,\varphi) \,{\to}\, (Q',\varphi')$ of $P$-twisted bundles}
is a morphism $f\colon Q \,{\to}\, Q'$ of $G_2$-bundles such that 
$\varphi' \,{\circ}\, f \,{=}\, \varphi$.
\\[3pt]
(iii)\,
The category of $P$-twisted $G_1$-bundles is denoted by $\mathcal{A}_{G_1}\big(P{\to}M\big)$.
\end{definition}

A $J$-equivariant topological field theory can now be obtained explicitly by a
direct generalization of the construction of ordinary Dijkgraaf-Witten theories.
In particular, for any group element $j \,{\in}\, J$
there is a $J$-cover $P_j$ of $\mathbb{S}^1$ with monodromy $j$ and thus a category 
$\mathcal{C}_j \,{=}\, [\mathcal{A}_G\big(P_j{\to}\mathbb{S}^1\big),\mathrm{vect}]$.
Then the category $\mathcal{C} \,{:=}\, \bigoplus_{j\in J}\mathcal{C}_j$
can be shown \cite{cs-mns} to be a $J$-equivariant modular category. Its orbifold 
category is the category of modules over the Drinfeld double of the group $G_2$.

The idea to extend the Dijkgraaf-Witten construction to categories
of generalizations of bundles turns out to be quite fruitful.
As an illustration consider topological field theories with defects. Such theories
are defined on cobordism categories consisting of manifolds with singularities,
see e.g.\ \cite[Sect.\,4.3]{cs-Lurie}. They turn out to be most intimately
linked to the theory of module categories over fusion categories
\cite{cs-FSV1}. Dijkgraaf-Witten theories with defects can be
realized using a generalization of relative bundles \cite{cs-FSV2}.
In this way a combination of field-theoretic and geometric arguments allows one 
to recover pure algebraic results \cite{cs-ostrik} on module categories over the 
fusion category of finite-dimensional vector spaces that are graded by a finite group.


\section{Partial dualizations in the theory of Hopf algebras}

The program, initiated by Andruskiewitsch and Schneider, of classifying 
finite-di\-men\-sional pointed Hopf algebras over a field has grown into one
of the more remarkable recent developments in algebra; see \cite{cs-Andr} for a
review. One of the surprises that were encountered in the 
realization of this program is a close connection between pointed Hopf algebras and 
Lie theory. The crucial ingredient of this connection is the theory of Weyl groupoids 
(see \cite{cs-hecke}, as well as \cite{cs-cuLe}).

In this groupoid a construction is central that, while being purely algebraic, turns out 
to have important ramifications in topological field theory. 
In the form introduced originally \cite{cs-heSc3}, this construction makes
extensive use of smash products; a convenient formulation
in the framework of braided monoidal categories proceeds as follows.

We first summarize a few facts that are needed in the construction \cite{cs-BLS}:
  \begin{itemize}
  \item
The notions of a Hopf algebra $A$, of a Hopf pairing 
$\omega\colon A\,{\otimes}\, B \,{\to}\, {\mathbf 1}_{\mathcal C}$, and of the 
category $_A^A \mathcal{YD}(\mathcal{C})$ of Yetter-Drinfeld modules over $A$ make sense
in any braided category $\mathcal C$. (Recall that a Yetter-Drinfeld module has the 
structure of both a module and of a comodule, subject to a certain 
consistency condition. Yetter-Drinfeld modules form a braided monoidal category.)

  \item
Suppose that we are given  two Hopf algebras $A$ and $B$ in $\mathcal C$ together with a
non-de\-ge\-ne\-rate Hopf pairing 
$\omega\colon A \,{\otimes}\, B \,{\to}\, {\mathbf 1}_{\mathcal C}$.
Using $\omega$ and its inverse, we can turn a left $A$-action into a left $B$-coaction and 
a left $A$-coaction into a left $B$-action. Hereby we obtain an isomorphism 
  $$
  \Omega^\omega:\quad{}_A^A \mathcal{YD}(\mathcal{C}) \,\stackrel\cong\longrightarrow\,
  {}_B^B \mathcal{YD}(\mathcal{C})
  $$
of braided categories.
\end{itemize}

The input of the construction performed in \cite{cs-BLS} is a {\em partial dualization datum}
$\mathcal A$ for a Hopf algebra $H$ in a braided category $\mathcal C$; this consists of 
a projection $\pi\colon H \,{\to}\, A$ 
to a Hopf subalgebra, and a Hopf algebra $B$ together with a non-degenerate
Hopf pairing $\omega\colon A \,{\otimes}\, B \,{\to}\, \mathbf{1}_{\mathcal C}$.
A partial dualization datum $\mathcal A$ gives rise to a new Hopf algebra 
$r_{\!\mathcal A}(H)$ in $\mathcal C$ as follows.

\begin{enumerate}
\item The Radford projection theorem, applied to the projection 
$\pi\colon H \,{\to}\, A$, allows us to write the Hopf algebra
$H$ as a Radford biproduct $H \,{\cong}\, K {\rtimes} A$, with $K$ the Hopf
algebra in the braided category $_A^A \mathcal{YD}(\mathcal{C})$ given
by the coinvariants of $H$ with respect to the projection $\pi$.

\item 
The image $L\,{:=}\, \Omega^\omega(K)$ of the Hopf algebra $K$ in the braided category 
$_A^A \mathcal{YD}(\mathcal{C})$ under the braided monoidal equivalence $\Omega^\omega$ 
is a Hopf algebra in the braided category $_B^B \mathcal{YD}(\mathcal{C})$.

\item 
The partially dualized Hopf algebra $r_{\!\mathcal{A}}(H)$ is defined as the Radford biproduct 
  $$
  r_{\!\mathcal A}(H) := L\rtimes B
  $$
of $L$. This is a Hopf algebra in the braided category $\mathcal{C}$.
\end{enumerate}

\noindent
In short, we simultaneously dualize a Hopf subalgebra $A$ of $H$ and 
transport the
coinvariants $K \,{\subseteq}\, H$ covariantly to $L\,{\subseteq}\, r_{\!\mathcal A}(H)$.
Note that this combines a contra- and a covariant operation; thus 
partial dualization is not functorial in the Hopf algebra $H$.

\medskip

{}From the partial dualization datum $\mathcal A$ for $H$ one can
obtain a partial dualization datum $\mathcal{A}^{\!-}_{\phantom|}$ 
of the Hopf algebra $r_{\!\mathcal{A}}(H)$. There is then a canonical isomorphism 
  $$
  r_{\!\!\mathcal{A}^{\!-\!}_{\phantom|}}\big( r_{\!\mathcal A}^{}(H) \big) \,\cong\, H 
  $$
of Hopf algebras in $\mathcal{C}$,
showing that partial dualization is essentially involutive.

The partially dualized Hopf algebra $r_{\!\mathcal{A}}(H)$ is rather different from
$H$, in general. As an illustration consider the Taft algebra $T_\zeta$ for some
primitive $d\,$th root of unity $\zeta$. $T_\zeta$ is a Hopf algebra over the complex
numbers that is generated by a group-like element $g$ of order $d$ and a skew-pri\-mi\-tive 
element $x$ having coproduct $\Delta(x) \,{=}\, g \,{\otimes}\, x \,{+}\, x \,{\otimes}\, 1$.
There is a projection $\pi$ from the Taft algebra $T_\zeta$ to the Hopf subalgebra 
$A \,{\cong}\, \mathbb{C}[\mathbb{Z}_d]$ generated by $g$. When taking this projection 
to obtain a partial dualization datum, the resulting partial dualization $r_{\!\mathcal A}(H)$
is in fact isomorphic to the Taft algebra itself, with the isomorphism depending on a 
choice of a Hopf pairing $\omega\colon A \,{\otimes}\, A \,{\to}\,\mathbb{C}$, and 
thereby on a choice of a primitive $d\,$th root of unity. 
By considering central extensions $\hat{T}_{\zeta,q}$ of $T_\zeta$ by group-like elements
one can get Hopf algebras with non-trivial partial dualization: the partial dualization 
$r_{\!\mathcal A}(\hat{T}_{\zeta,q})$ no longer contains any central group-like elements.
Rather, the coproduct of the skew-pri\-mi\-tive element of 
$\check{T}_{\zeta,q}$ gets modified, in such a way its partial dualization
has additional central characters.

A priori the categories of left modules in $\mathcal C$ over
the Hopf algebra $H$ and over its partial dualization $r_{\!\mathcal A}(H)$ 
are thus rather different. However, there does exist a representation-theoretic relation
between them \cite{cs-BLS}: The categories of Yetter-Drin\-feld modules over a Hopf algebra
$H$ in $\mathcal{C}$ and of those over $r_{\!\mathcal A}(H)$ are braided equivalent: we have
  $$
  {}_H^H \mathcal{YD} (\mathcal{C})
  \,\cong\, {}_K^K \mathcal{YD}\big( {}_A^A \mathcal{YD}(\mathcal{C}) \big)
  \,\stackrel{\Omega^\omega}{\longrightarrow}\,
  {}_L^L\mathcal{YD}\big( {}_B^B \mathcal{YD}(\mathcal{C}) \big)
  \,\cong\, {}_{r_{\!\mathcal A}(H)}^{r_{\!\mathcal A}(H)}\mathcal{YD}(\mathcal{C}) \,.
  $$
Here the first and third isomorphisms come from the isomorphism
  $$
  {}_{K{\rtimes}A}^{K{\rtimes}A}\mathcal{YD}(\mathcal{C})
  \,\cong\, {}_K^K \mathcal{YD} \big( {}_A^A \mathcal{YD}(\mathcal{C}) \big)
  $$
of braided categories for a biproduct.

The equivalence of categories of Yetter-Drinfeld modules implies a relation between 
the Hopf algebra $H$ and and its partial dualization $r_{\!\mathcal A}(H)$. Let us discuss
this for the case of a finite-dimensional Hopf algebra over a field. The category of 
Yetter-Drinfeld modules over $H$ is then the Drinfeld center of the category of $H$-modules. 
Now semisimple algebras having isomorphic centers are Morita equivalent. This statement 
allows for a categorification: the notion of an algebra is substituted with the one 
of a monoidal category, while the one of a module over an algebra gets replaced by
the one of a module category. This way one arrives \cite[Thm\,3.1]{cs-ENO11} at the statement
that semisimple tensor categories with braided-equi\-va\-lent Drinfeld centers 
have equivalent bicategories of module categories, i.e.\ they are related by what is called
\cite{cs-Mue03} a weak monoidal Morita equivalence. It is natural to ask whether 
the bicategories of module categories over $H$-mod and over $r_{\!\mathcal A}(H)$-mod 
are closely related, possibly even equivalent. This in turn suggests that the 
categorified representation theory of a Hopf algebra $H$, i.e.\
the theory of module and bimodule categories over the monoidal category
$H$-mod, should be regarded as an important aspect of the theory of Hopf algebras.

Since the Brauer-Picard group of the monoidal category $H$-mod is isomorphic to the 
group of braided equivalences of the Drinfeld center $\mathcal Z(H\textrm{-mod})$
\cite{cs-enoM}, it is invariant under partial dualization. On the other hand  it is
not hard to give conditions \cite{cs-LentPriel} that ensure the existence of a Hopf 
isomorphism $f\colon r_{\!\mathcal A}(H) \,{\stackrel\cong\to}\, H$. 
(There is also a generalization involving a modified partial dualization of $H$.)
Every pair $(r_{\!\mathcal A}, f)$ then induces an equivalence
  $$
  \mathcal Z(H\textrm{-mod}) \,\simeq\, \mathcal Z(r_{\!\mathcal A}(H)\textrm{-mod})
  \,\simeq_f\, \mathcal Z(H\textrm{-mod}) 
  $$
and thus provides an element in the group of braided autoequivalences of the Drinfeld center
$\mathcal Z(H\textrm{-mod})$, and thereby an element of the Brauer-Picard group.

This result has an interesting application in three-dimensional
topological field theories of Turaev-Viro type. A Turaev-Viro theory is a
fully extended oriented three-di\-men\-sional topological field theory,
which assigns a spherical fusion category $\mathcal A$ to the point.
Types of topological surface defects in this theory are given by
$\mathcal A$-bi\-mo\-dule categories. Now it is a general field-theoretic
pattern that topological invertible codimension-one defects are related to
symmetries. As a consequence, the Brauer-Pi\-card groups of the
fusion categories in question
can be identified \cite{cs-FPSV} as symmetry groups of three-di\-men\-sional topological
field theories of Turaev-Viro type. In the special case that the fusion category 
is the category of modules over the group algebra $H \,{=}\, \mathbb{C}[G]$ of 
an abelian group $G$, the reflections just discussed amount to what in the
physics literature are known as electric-magnetic dualities. This constitutes yet another
instance of an enlightening relation between algebraic notions and notions from field theory.

\vfill

\pagebreak 
  

\pagebreak   

\section*{References}

\renewcommand{\refname}{}    
\vspace*{-26pt}              

\end{document}